\documentclass{ifacconf}

\usepackage{graphicx}      % include this line if your document contains figures
\usepackage{natbib}        % required for bibliography
\usepackage{amsmath,amssymb}
\usepackage{tikz,pgf}
\usetikzlibrary{decorations.markings, decorations.pathmorphing, positioning}
\tikzstyle{every node}=[circle, draw,inner sep=0pt, minimum width=4pt]

\newcommand{\skp}[1]{\big\langle #1 \big\rangle}

\newcommand{\funSpace}[1]{\mathcal{#1}} 
\newcommand{\Vsa}{\funSpace{Q}}
\newcommand{\Vsb}{\funSpace{V}}
\begin{document}
\begin{frontmatter}

\title{Certified Reduced Basis Method for the Damped Wave Equations on Networks}% \thanksref{footnoteinfo}} 

%\thanks[footnoteinfo]{Sponsor and financial support acknowledgment
%goes here. Paper titles should be written in uppercase and lowercase
%letters, not all uppercase.}

\author[UT]{Nadine Stahl} 
\author[UT]{Bj{\"o}rn Liljegren-Sailer}
\author[UT]{Nicole Marheineke} 

\address[UT]{Trier University, 
   54286 Trier, Germany (e-mail: \{nadine.stahl, bjoern.sailer, marheineke\}@uni-trier.de).}

\begin{abstract}                % Abstract of not more than 250 words.
%These instructions give you guidelines for preparing papers for the $10^{\mathrm{th}}$ Vienna Conference on Mathematical Modelling. Please use this document as a template to prepare your manuscript. For submission guidelines, follow instructions on paper submission system as well as the event website. 
In this paper we present a reduced basis method which yields structure-preservation and a tight a posteriori error bound for the simulation of the damped wave equations on networks. The error bound is based on the exponential decay of the energy inside the system and therefore allows for sharp bounds without the need of regularization parameters. The fast convergence of the reduced solution to the truth solution as well as the tightness of the error bound are verified numerically using an academic network as example.
\end{abstract}


\end{frontmatter}

\section{Introduction}

The dynamics of many technical applications can be described by partial differential equations. Depending on the complexity of the application, a high fidelity approximation may require a very fine discretization and thus leads to high-dimensional systems of equations.
The computational burden gets very high in these cases, which can limit the applicability, especially for many-query tasks in optimization and control. As a remedy, model order reduction methods have been developed. The reduced basis method is a model order reduction approach specialized for parameter-dependent settings. Using a greedy procedure, reduced low-order models of a user-pre-defined accuracy are constructed. While training of the reduced models might be costly, online evaluations performed repeatedly during many tasks are cheap and very efficient. The crucial ingredient of the reduced basis approach is the derivation of a problem-adapted, tight a posteriori error bound.

\cite{Grepl2005} derived respective bounds for parabolic equations, which \cite{Haasdonk2008} used to set up efficient offline-online decomposed reduced basis methods. An early work on reduced basis methods in fluid dynamics is \citep{Veroy2003}, where the Burger's equation is considered. \cite{Knezevic2010} developed a similar strategy for the Boussinesq equations. The Stokes equations in a two-dimensional setting are treated in \citep{Gerner2012, Gerner2012a}, where the ideas of the offline-online splitting have also been adopted.

In this paper we deal with the damped wave equations on a network of edges. This model problem can be used to describe acoustic waves or a simplified gas pipeline network. While high-order approximations on a single edge can be constructed by standard discretization methods, the network aspect makes the problem particularly interesting and challenging for model order reduction. Similar equations have been considered in the reduced basis context, however, not on networks.
Works on the undamped wave equations as well as others on Langrangian hydrodynamics can be found in, e.g., \citep{Amsallem2014, Glas2020, Copeland2021}. Existing approaches for hyperbolic reduced order models with error bound calculate approximative inf-sup-parameters or use penalty terms, cf., \citep{Gerner2012, Gerner2012a}. Moreover, generic approaches for ordinary differential equations exist, see, e.g., \citep{Haasdonk2011}. 

Our approach strongly relies on the analytical results from \citep{Egger2018a, Egger2019, Kugler2019}. The papers show exponential stability for the damped wave equations on networks and certain Galerkin approximations of these equations. The derivation of our a posteriori error bound is based on the latter results. 
While \cite{Egger2018} developed a structure-preserving model order reduction based on Krylov methods for the linear damped wave equations, \cite{Sailer2021, Sailer2021a} used a structure-preserving proper orthogonal decomposition and extended the model problem to nonlinear flows on networks. In both cases, structure preservation during the reduction process was achieved by constructing compatible spaces, yielding stable and more robust reduced models.
In this paper we aim to transfer and embed the proposed structure-preserving reduction strategy into the reduced basis method. In the reduced basis method, a hierarchy of Galerkin approximations is employed. The highest order, the so-called truth solution, is obtained in our approach by a mixed finite element method. The truth is approximated by lower-dimensional and cheaper-to-evaluate reduced basis solutions. The novelty of our approach, which distinguishes it from conventional reduced basis methods, is that it takes into account compatibility constraints during training and thus can guarantee structure-preserving and stable reduced models. It also allows for a problem-adapted error bound that performs better than the more generic error bound by \cite{Haasdonk2011}. We illustrate the fast convergence of the reduced order model and the tightness of our error bound by numerical results.

%%%%%%%%%%%%%%
%%%%%%%%%%%%%%
\section{Problem Setting}

A network of pipes is modeled by a directed graph $\mathcal{G}(\mathfrak{V}, \mathfrak{E})$  with each edge $e \in \mathfrak{E}$ representing a pipe of length $l^e>0$. The set of nodes $\mathfrak{V}$ divides into junctions and boundary nodes, where the latter are the nodes only incident to one edge.
Moreover, we define the sets of all topologically ingoing and outgoing edges to a node $v$ by
\begin{equation*}
\begin{split}
\delta_v^+ &= \{e \in \mathfrak{E}: \exists w \in \mathfrak{V} \text{ with } e = (v,w)\},\\
\delta_v^- &= \{e \in \mathfrak{E}: \exists w \in \mathfrak{V} \text{ with } e = (w,v)\}.
\end{split}
\end{equation*}
To describe the system dynamics, we identify the edges $e \in \mathfrak{E}$ with intervals $[0,l^e]$. The edgewise pressure $p^e(t,x)$ and mass flux $u^e(t,x)$ for $(t,x) \in [0,T] \times [0,l^e]$ are governed by the damped wave equations
\begin{equation} \label{eq:dwe_strong}
\begin{split}
a^e\partial_t p^e + \partial_x u^e &= f^e,\\
b^e\partial_t u^e + \partial_x p^e + d^e u^e &= g^e.
\end{split}
\end{equation}
The pipe constants  $a^e,b^e, d^e$ are allowed to vary over a parameter domain. Our system is thus parametrized in these constants, although we do not indicate this in the notation to keep it more concise. 

The pipe equations are coupled at the junctions via the Kirchhoff conditions
\begin{eqnarray} 
&\sum_{e \in \delta_v^-} u^e(l^e,t) = \sum_{e \in \delta_v^+} u^e(0,t),\label{eq:kirchflux}\\
&p^e(l^e,t) = p^{v}(t), e \in \delta_v^+, \quad  p^e(0,t) = p^{v}(t), e \in \delta_v^-, \nonumber
\end{eqnarray}
with the auxiliary variable $p^{v}$ denoting the pressure at node $v$. To close the network model, we assume homogeneous pressure boundary conditions and consistent initial values $p_0(x), u_0(x)$ to be prescribed. 
For ease of presentation, we only treat homogeneous boundary conditions in this paper. Other boundary conditions could be treated similarly, but note that, by construction, only a system with homogeneous boundary conditions needs to be considered in the derivation of our a posteriori error bound (Theorem~\ref{thm:rbbound}) anyway.

We make the following assumption for the model \eqref{eq:dwe_strong}.
\begin{assum} \label{assum:constants}
There exist constants $C_0, C_1 > 0$, such that $C_0 < a^e,b^e,d^e < C_1$ for $e \in \mathfrak{E}$. The right hand sides $f^e$, $g^e$ satisfy $f^e, g^e \in \mathcal{L}^{1}(0,T; \mathcal{L}^2(e))$ and $\partial_t f^e, \partial_t g^e \in \mathcal{L}^{1}(0,T; \mathcal{L}^2(e))$ for $e \in \mathfrak{E}$.
\end{assum}

We consider our model problem on network generalizations of standard functions spaces. Given the Lebesgue space $\mathcal{L}^2(e)$  of square integrable functions on $e \in \mathcal{E}$, we define 
\begin{equation*}
\begin{split}
\mathcal{L}^2(\mathfrak{E}) &= \{u:\mathfrak{E}\rightarrow \mathbb{R},\, u|_e \in \mathcal{L}^2(e), \, e \in \mathfrak{E}\},\\
\mathcal{H}^{div}(\mathfrak{E}) &= \{u \in \mathcal{L}^2(\mathfrak{E}): \partial_x u|_e \in \mathcal{L}^2(e)  \text{ and } \eqref{eq:kirchflux} \text{ holds}\}.
\end{split}
\end{equation*}
The functions $u \in \mathcal{H}^{div}(\mathfrak{E})$ are only piecewise continuous. Thus, when we write $\partial_x u$, we refer to the broken (edgewise) derivative. The function evaluations in \eqref{eq:kirchflux} are well-defined in the sense of the trace theorem. The $\mathcal{L}^2$-scalar product on the network and its induced norm are given by
\begin{equation*}
\skp{f,g} =\sum_{e \in \mathfrak{E}} \int_e f^e g^e dx, \qquad \|\cdot\| = \skp{\cdot,\cdot}^{1/2}.
\end{equation*}
Here and in the following, network generalizations of functions $ f^e$ or constants $a^e,b^e,d^e$ are denoted with the same symbol but without the super-index $e$. Moreover, for a separable reflexive Banach space $\mathcal{X}$ we consider the Sobolev space 
\begin{equation*}
\mathcal{W}(0,T;\mathcal{X}) = \{v \in \mathcal{L}^2(0,T; \mathcal{X}): \partial_t v \in \mathcal{L}^2(0,T; \mathcal{X}')\},
\end{equation*}
where $\mathcal{X}'$ denotes the dual space of $\mathcal{X}$.

We can now state the variational formulation of our network problem.
\begin{prob}[Variational Formulation]\label{prob:con} Find $(p,u) \in \linebreak \mathcal{W}(0,T; \mathcal{L}^2(\mathfrak{E}) \times \mathcal{H}^{div}(\mathfrak{E}))$ with $p(0) = p_0, u(0) = u_0$ and such that for all $(q,v) \in \mathcal{L}^2(\mathfrak{E}) \times \mathcal{H}^{div}(\mathfrak{E})$ and a.e. $t \in (0,T]$
\begin{equation*}
\begin{split}
\skp{a\partial_t p, q} + \skp{\partial_x u, q} &= \skp{f,q}, \\
\skp{b\partial_t u, v} - \skp{p, \partial_x v} + \skp{du,v} &= \skp{g,v}.
\end{split}
\end{equation*}
\end{prob}
Problem \ref{prob:con} is well-posed and has a unique solution, see \citep{Egger2018a}.

%%%%%%%%%%%%%%%%%%
%%%%%%%%%%%%%%%%%%
\section{Hierarchy of approximations}

The reduced basis approach relies on a hierarchy of Galerkin approximations for Problem~\ref{prob:con}. Its highest fidelity model is the 'truth' solution given by a finite element discretization. As this model is typically too high-dimensional and computationally too expensive for many-query tasks, we consider lower-order and cheaper-to-evaluate approximations given by reduced basis projections. We particularly choose the truth solution and the reduced basis approximations such that the results from \cite{Egger2018a} apply to them, which makes the approach structure-preserving.

For the truth approximation, we assume a partitioning 
of every $e \in \mathfrak{E}$ to be given and denote by $\mathcal{P}^k(e)$ the space of piecewise polynomial functions of degree $k$ on that partitioning. The truth Galerkin spaces are given by
\begin{equation*}
\begin{split}
\mathcal{Q} &= \{q \in \mathcal{L}^2(\mathfrak{E}), \, q_{|e} \in \mathcal{P}^0(e), \, e \in \mathfrak{E}\},\\
\mathcal{V} &= \{v \in \mathcal{H}^{div}(\mathfrak{E}), \, v_{|e} \in \mathcal{P}^1(e), \, e \in \mathfrak{E}\}.
\end{split}
\end{equation*}
The functions in $\mathcal{Q}$ are piecewise constant, whereas $\mathcal{V}$ consists of piecewise linear functions that are continuous on every edge. At the junctions the Kirchhoff conditions \eqref{eq:kirchflux} hold.
\begin{prob}[Truth Solution] \label{prob:truth} Find $(p,u) \in \mathcal{W}(0,T; \mathcal{Q} \times \mathcal{V})$ with $p(0) = p_0\in \mathcal{Q}$, $u(0) = u_0\in \mathcal{V} $ and 
\begin{equation*}
\begin{split}
\skp{a \partial_t p, q} + \skp{\partial_x u, q} &= \skp{f,q}, \\
\skp{b \partial_t u, v} - \skp{p, \partial_x v} + \skp{du,v} &= \skp{g,v},
\end{split}
\end{equation*}
for all $(q,v) \in \mathcal{Q}\times \mathcal{V}$ and a.e. $t \in (0,T]$.
\end{prob}
The following a priori bound can be shown for the truth solution.
\begin{thm}[\cite{Egger2018a}]\label{thm:egger}
Let $(p,u)$ be a solution of Problem~\ref{prob:truth}. Then for $t \geq s \geq 0$ 
\begin{equation*}
\begin{split}
\|p(t)\|^2 + \|&u(t)\|^2 \leq C' \exp({-\gamma(t-s)})(\|p(s)\|^2+\|u(s)\|^2) \\
&+ C'' \int_s^t \exp({-\gamma(t-r)})(\|f(r)\|^2+\|g(r)\|^2)dr
\end{split}
\end{equation*}
with stability constants $\gamma$, $C'$, $C''$ independent of $s$, $t$ and of the data $f$, $g$.
\end{thm}
Note that the constants $\gamma$, $C'$, $C''$ can be explicitly determined, which will be crucial for our approach. Theorem~\ref{thm:egger} implies the truth solution to be well-posed and exponentially stable. The same can be shown for more general Galerkin approximations that satisfy certain compatibility conditions, see \citep{Egger2018a,Egger2018}. This motivates the development of reduced basis approximations that are compatible.

A reduced basis approximation is established by finding suitable subspaces $\Vsa_N \subset \Vsa$ and $\Vsb_N \subset \Vsb$, which are of much lower dimension, $N = \dim(\Vsa_N \times \Vsb_N) \ll \dim(\Vsa \times \Vsb)$. The aforementioned compatibility conditions on the ansatz spaces are now presented. Note that they are naturally fulfilled by the truth solution spaces $\Vsa$ and $\Vsb$. 

\begin{assum}[Compatibility Conditions]\label{assum:compSpaces}
The an\-satz \linebreak spa\-ces $\Vsa_N \subset \Vsa$ and $\Vsb_N \subset \Vsb$ are such that
\begin{enumerate} 
\item[A1)] $\Vsa_N =    \left \{\xi: \text{It exists } \zeta \in \Vsb_N \text{ with } \partial_x \zeta = \xi \right \}$,
\item[A2)] $\funSpace{K} \subset \Vsb_N$, \hspace{0.14cm} $\funSpace{K}=\{v \in \funSpace{H}^{div}(\mathfrak{E}): \, \partial_x v = 0 \}$.
\end{enumerate}
\end{assum}

Note that the space $\funSpace{K}$ is low-dimensional, as it consists of edgewise constant fluxes. We are now able to formulate a well-posed and stable reduced basis approximation.
\begin{prob}[Reduced Basis Approximation] \label{prob:rb} Let $\Vsa_N \times \Vsb_N$ fulfill Assumption~\ref{assum:compSpaces}. Find $(p^N,u^N) \in \mathcal{W}(0,T; \Vsa_N \times \Vsb_N)$ with $p^N(0) =  \Pi_{\Vsa_N} p_0$, $u^N(0) = \Pi_{\Vsb_N}u_0$ and
\begin{equation*}
\begin{split}
\skp{a \partial_t p^N, q^N} + \skp{\partial_x u^N, q^N} &= \skp{f,q^N}, \\
\skp{b \partial_t u^N, v^N} - \skp{p^N, \partial_x v^N} + \skp{du^N,v^N} &= \skp{g,v^N},
\end{split}
\end{equation*}
for all $(q^N,v^N) \in \Vsa_N\times \Vsb_N$ and a.e. $t \in (0,T]$, with $\Pi_{\Vsa_N}$, $\Pi_{\Vsb_N}$ denoting $\funSpace{L}^2$-projections onto the reduced spaces.
\end{prob}

%%%%%%%%%%%%%%%%%%%
%%%%%%%%%%%%%%%%%%%
\section{Reduced Basis Approximation}

The reduced basis approach consists of the reduced basis approximation itself and an error bound, which both need to be evaluable by an efficient offline-online decomposition. In the offline-phase, first the parameter- and solution-independent parts of the error bound are prepared, and then a reduced model is set up by training it towards the parameter range of interest. The latter step relies on a greedy procedure employing the error bound. Note that the reduced models inherit the parameter-dependence of the truth solution, but the reduction bases are constructed to be parameter-independent. The unique and novel feature of our approach, which distinguishes it from the conventional reduced basis methods, is that it regards the compatibility conditions in the training and thus can guarantee structure-preserving and stable reduced models.

%%%%%%
\subsection{A Posteriori Error Bound}

Let $(p,u)$ and $(p^N, u^N)$ be solutions to the truth solution and the reduced basis approximation, respectively.
The error we like to control reads
\begin{equation}\label{eq:def_err}
\begin{split}
e^p &= p - p^N \in \mathcal{L}^2(0,T;\mathcal{Q}),\\
e^u &= u - u^N \in \mathcal{L}^2(0,T;\mathcal{V}).
\end{split}
\end{equation}
We derive a residual-based error bound with the residual $(r^{p}(t),r^{u}(t)) \in \mathcal{Q} \times \mathcal{V}$, $t \geq 0$, defined by
\begin{equation}\label{eq:def_res}
\begin{split}
\skp{r^{p}, q} &= \skp{f,q} - \skp{a \partial_t p^N, q} - \skp{\partial_x u^N, q}, \\
\skp{r^{u}, v} &= \skp{g,v} - \skp{b \partial_t u^N, v} + \skp{p^N_, \partial_x v} -\skp{du^N,v}
\end{split}
\end{equation}
for $(q,v) \in \mathcal{Q} \times \mathcal{V}$. Moreover, we make use of a generalized Poincare constant $C_P$, defined as the optimal constant fulfilling
\begin{equation} \label{eq:poincare}
\|b^{1/2}u\|^2 \leq C_P^2(\|a^{-1/2}\partial_x u\|^2+\|d^{1/2}\Pi_0 u\|^2)
\end{equation}
for all $u \in \Vsb$, with $\Pi_0:\Vsb \rightarrow \funSpace{K}$ given as the $\funSpace{L}^2$-projection onto the space of constant fluxes $\funSpace{K}$ from Assumption~\ref{assum:compSpaces}.

\begin{thm}[A-posteriori Error Bound] \label{thm:rbbound}
Let $(p,u) \in \break \mathcal{W}(0,T; \mathcal{Q} \times \mathcal{V})$ and $(p^N,u^N) \in \mathcal{W}(0,T;\mathcal{Q}_N \times \mathcal{V}_N)$ be solutions to Problem~\ref{prob:truth} and Problem~\ref{prob:rb}, respectively. Then the error~\eqref{eq:def_err} fulfills
\begin{equation*}
\begin{split}
\|e^p(t)\|^2 + \|e^u(t)\|^2 \leq \Delta(t) , \quad \text{ for $t \geq 0$}, 
\end{split}
\end{equation*}
with
\begin{equation*}
\begin{split}
\Delta(t) &= C' \exp({-\gamma t})(\|e^p(0)\|^2+\|e^u(0)\|^2) \\
&\quad + C'' \int_0^t \exp({-\gamma(t-\tau)})(\|r^{p}(\tau)\|^2+\|r^{u}(\tau)\|^2)d\tau.
\end{split}
\end{equation*}
The stability constants can be estimated in terms of the constants defined in Assumption~\ref{assum:constants} and \eqref{eq:poincare} by
\begin{equation*}
\begin{split}
C'=C''=3\left(\frac{C_1}{C_0}\right)^{1/2},  \quad \gamma = \frac{2}{3} \frac{C_0}{C_1}\frac{C_0}{2C_0+4C_PC_1}.
\end{split}
\end{equation*}
\end{thm}
\begin{pf} 
The residual equation \eqref{eq:def_res} can be written as a differential equation in the error $e^p$ and $e^u$, i.e., for $q \in \mathcal{Q}$ and $ v \in \mathcal{V}$,
\begin{equation*} %\label{eq:error_pde}
\begin{split}
\skp{a \partial_t e^p, q} + \skp{\partial_x e^u, q} &= \skp{r^{p}, q},\\
\skp{b \partial_t e^u, v} - \skp{e^p, \partial_x v} + \skp{de^u,v} &= \skp{r^{u}, v}.
\end{split}
\end{equation*}
This shows that $e^p$ and $e^u$ solve Problem~\ref{prob:truth}  with right hand side functions $r^{p}$ and $r^{u}$. Hence, we can apply Theorem~\ref{thm:egger}. The explicit expressions for the stability constants are shown in \citep{Egger2018}. \qed
\end{pf}
Let us emphasize the online-efficiency of our error bound, which holds despite the parameter dependence of our model problem. Under the usual assumption that the parameter dependence (in the constants $a,b,d$) allows for an affine-linear representation, the residual norms $\|r^{p}\|$ and $\|r^{u}\|$ can be evaluated using an efficient offline-online decomposition, cf., \citep{Grepl2005}. The stability constants $C',C''$ and $\gamma$ are determined once in the offline phase. The computationally most demanding ingredient is the Poincare constant $C_P$ that is needed for the evaluation of $\gamma$. It is given as the largest eigenvalue of the generalized (and parameter-dependent) eigenvalue problem
\begin{align}\label{eq:cp_eig}
Bu = \lambda (A+D)u,
\end{align}
cf., \citep{Egger2018}, with $B,A$ and $D$ denoting the matrix representations of the operators $\tilde{B},\tilde{A}$ and $\tilde{D}$, which, in turn, are defined by
\begin{equation*}
\begin{split}
\skp{\tilde{B}u,v} &= \skp{bu,v}, \qquad \skp{\tilde{A}u,v}=  \skp{a^{-1}\partial_x u, \partial_x v},\\
\skp{\tilde{D} u,v} &= \skp{d \Pi_0 u, \Pi_0 v}
\end{split}
\end{equation*}
for $u,v \in \Vsb$ and $\Pi_0$ as in \eqref{eq:poincare}.

\begin{rem}
As an alternative to solving the eigenvalue problem~\eqref{eq:cp_eig}, $\gamma$ can also be estimated by a simulation-based approach, cf., \citep{Egger2018a}. As shown there, solutions to compatible Galerkin approximations with constant right hand sides and homogeneous boundary conditions satisfy
\begin{equation*}
\mathcal{E}(t) \leq C' \exp(-\gamma(t-s))\mathcal{E}(s), \quad \text{for }t\geq s,
\end{equation*}
with $\mathcal{E}(t) = 0.5(\|a^{1/2}\partial_t p(t)\|^2 + \|b^{1/2}\partial_t u(t)\|^2)$ and the constants as in Theorem~\ref{thm:rbbound}. Given a solution trajectory, $\gamma$  is determined by a least squares fit to the values of $\mathcal{E}(t)$ for $t \geq s$ with fixed $s$. For details on the approach, we refer to the reference. Note that in practice, good estimates can already be obtained from simulations of our reduced models, as they inherit similar stability constants due to the compatible Galerkin ansatz.
\end{rem}

%%%%%%
\subsection{Compatible Reduced Basis Spaces}

Our compatible reduced basis approximations are constructed with a POD greedy algorithm similar to \citep{Haasdonk2008}. However, in one point it fundamentally differs from the standard procedure. That is, we employ a constrained principal component analysis in order to fulfill Assumption~\ref{assum:compSpaces}. We briefly explain the procedure in the following. The details on the derivation and efficient algorithmic realization of the constrained principal component analysis can be found in \cite{Sailer2021a}. Let snapshots $(p_\ell,u_\ell) \in \Vsa \times \Vsb$, $\ell =1,\ldots, L$ of the truth solution (Problem~\ref{prob:truth}) to certain time points be given, i.e., $p_\ell = p(t_\ell)$, $u_\ell = u(t_\ell)$. A reduced space $\Vsa_N$  for the pressure is obtained by a principal component analysis of the joint snapshot collection  $\{p_1,\ldots p_L,\partial_x u_1, \ldots, \partial_x u_L \} \subset \Vsa$. Then, the compatible space $\Vsb_N$ for the mass flux is constructed from the compatibility conditions, i.e., $\Vsb_N = \partial_x^+\Vsa_N \oplus \funSpace{K}$ with $\partial_x^+$ denoting an arbitrary right-inverse of $\partial_x: \Vsb \rightarrow \Vsa$ and $\funSpace{K}$ as in Assumption~\ref{assum:compSpaces}. Note that the particular choice of the right-inverse has no influence on the approximation results. The resulting reduced space $\Vsa_N \times \Vsb_N$ of dimension $N$ is compatible and fulfills an optimality condition for the snapshot data under the compatibility constraint. The procedure has been used in \citep{Sailer2021a} for the barotropic Euler equations on networks and shown to give significantly better results than non-structure-preserving alternatives. 

To deal with the parameter-dependence of our problem, we embed the compatibility procedure in the greedy algorithm of \citep{Haasdonk2008}. Given a compatible reduced space $\Vsa_{N_i} \times \Vsb_{N_i}$, we compute the reduced solutions and the error bounds of Theorem~\ref{thm:rbbound} for all parameter settings in the training set. The parameter setting, for which the error bound indicates the worst approximation of the reduced solution to the truth, is selected to enrich the space. The temporal snapshots $(p_\ell,u_\ell)$ of the corresponding truth solution that do not already contribute to $\Vsa_{N_i}$ are used to build up $\tilde{\Vsa}$. The new compatible spaces are then $\Vsa_{N_{i+1}} = \tilde{\Vsa} \oplus \Vsa_{N_i}$ and $\Vsb_{N_{i+1}} = \partial_x^+\tilde{\Vsa} \oplus \Vsb_{N_i}$. We iterate until either a given maximum size of the reduced basis space is reached or the error bound meets a given tolerance for all training parameter settings.

\begin{rem}
Other approaches that handle equations for flow problems with a mixed formulation exist. In \citep{Gerner2012, Gerner2012a}, e.g., the required compatibility conditions on the reduced spaces are similar. But the reduced spaces are built by a separate principal component analysis for density and mass flux, respectively. Therefore, their approach does not lead to spaces that fulfill an optimality condition under the snapshot data. For details we refer the reader to \citep{Sailer2021a}.
\end{rem}

%%%%%%%%%%%%%%
%%%%%%%%%%%%%%
\section{Numerical Results}

\begin{figure}[tpb]
\centering
\tikzset{->-/.style={decoration={ markings, mark=at position 0.6 with {\arrow{>}}},postaction={decorate}}}
\begin{tikzpicture}[scale = 1]
	\draw[->-] (0,0) node[fill = black]{} -- (1.5,0) node[draw=none,fill=none, midway,sloped,above] {$e_1$};
	\draw[->-] (1.5,0)  -- (3,1.25) node[draw=none,fill=none, midway,sloped,above] {$e_2$};
	\draw[->-] (1.5,0) node[fill = black]{}-- (3,-1.25) node[draw=none,fill=none, midway,sloped,below] {$e_3$};
	\draw[->-] (3,1.25) -- (3,-1.25) node[draw=none,fill=none, midway,sloped,above] {$e_4$};
	\draw[->-] (3,1.25) node[fill = black]{}-- (4.5,0) node[draw=none,fill=none, midway,sloped,above] { $e_5$};
	\draw[->-] (3,-1.25) node[fill = black]{} -- (4.5,0) node[draw=none,fill=none, midway,sloped,below] { $e_6$};
	\draw[->-] (4.5,0) node[fill = black]{}-- (6,0) node[fill = black]{} node[draw=none,fill=none, midway,sloped,above] { $e_7$};
	\draw[draw=none] (0,-	0.1) node[below, draw=none, fill = none]{ $v_1$};
	\draw[draw=none] (6,-	0.1) node[below, draw=none, fill = none]{ $v_2$};
	%\draw (4,0.5) node[above , draw=none, fill = none]{\scriptsize $v_3$};
\end{tikzpicture}
\caption{Diamond network topology.}\label{fig:diamond}
\end{figure}
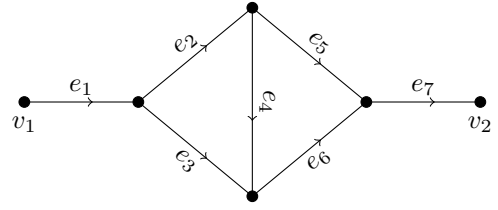%

The performance of our reduced order models and the a posteriori error bound is demonstrated for an academic network example, see diamond network in Figure \ref{fig:diamond}. The network consists of seven edges with length $l^e = 1$ for all $e \in \mathfrak{E}$. The edge parameters are chosen as
\begin{equation*}
\begin{split}
a &= \begin{bmatrix}4 &4& 1 &1 &1& 4 &4\end{bmatrix},\\ 
b &= \begin{bmatrix}0.25 &0.25& 1 &1& 1 &0.25 &0.25\end{bmatrix},\\ 
d &= \mu \begin{bmatrix}0.5& 0.5& 4 &4 &4& 0.5 &0.5\end{bmatrix},
\end{split}
\end{equation*}
where we, for ease of presentation, only consider one free parameter $\mu \in \mathfrak{P} = [0.01,10]$. 
The boundary conditions are taken as
\begin{equation*}
p^{v_1}(t) = 1- \cos(t),  \quad p^{v_2}(t) = 0,
\end{equation*}
and the initial conditions are set to zero. We consider the time interval $[0, 20]$ and perform the time integration of the truth and reduced models by means of the implicit L-stable Euler method with constant step size $\tau = 0.02$. This integrator is not structure-preserving, but dissipative. However, the resulting time-discretization error can be assumed to be negligibly small due to the chosen small step size. Our truth model is of size $1403$, which corresponds to $100$ finite elements per edge. For our error bound (Theorem \ref{thm:rbbound}) we compute the Poincare constant $C_P$ with the help of \eqref{eq:cp_eig}. Note that as $C_P$ and the model constants $C_0 = \min\{a^e,b^e,d^e\}$, $C_1 = \max\{a^e,b^e,d^e\}$ are dependent on the free parameter $\mu$, this also holds for the stability constants $C', C''$ and $\gamma$.

In training we use $12$ sample values for $\mu$ logarithmically distributed over the parameter domain $\mathfrak{P}$ to set up the compatible reduced basis space $\Vsa_N \times \Vsb_N$. The set of snapshots is given from the time discretization.
For testing the approximation quality, we consider a sample of 20 randomly chosen parameter values in $\mathfrak{P}$. For each parameter value we compute the actual error of the reduced approximation as well as our error bound $\Delta$.
For comparison, we consider the established, more generic residual-based error bound from \cite{Haasdonk2011}. It is of the form
\begin{equation*}
\begin{split}
\|e^p(t)\|^2 + \|e^u(t)\|^2 \leq \tilde{\Delta}(t) , \quad \text{ for $t \geq 0$}, 
\end{split}
\end{equation*}
with
\begin{equation*}
\begin{split}
\tilde{\Delta}(t) &= \tilde{C}\,\Big(\|e^p(0)\|+\|e^u(0)\| + \int_0^t \|r^{p}(\tau)\|+\|r^{u}(\tau)\| d\tau \Big)^2
\end{split}
\end{equation*}
whereby the constant $\tilde{C}$ is determined from the algebraic representation of truth and reduced solutions. Since our bound $\Delta$ inherits an exponentially decaying term under the time integral that $\tilde \Delta$ does not have, it is better adapted to the problem. Its effect can be well observed in the numerical results.

\begin{figure}[tpb]
\includegraphics[width = \columnwidth]{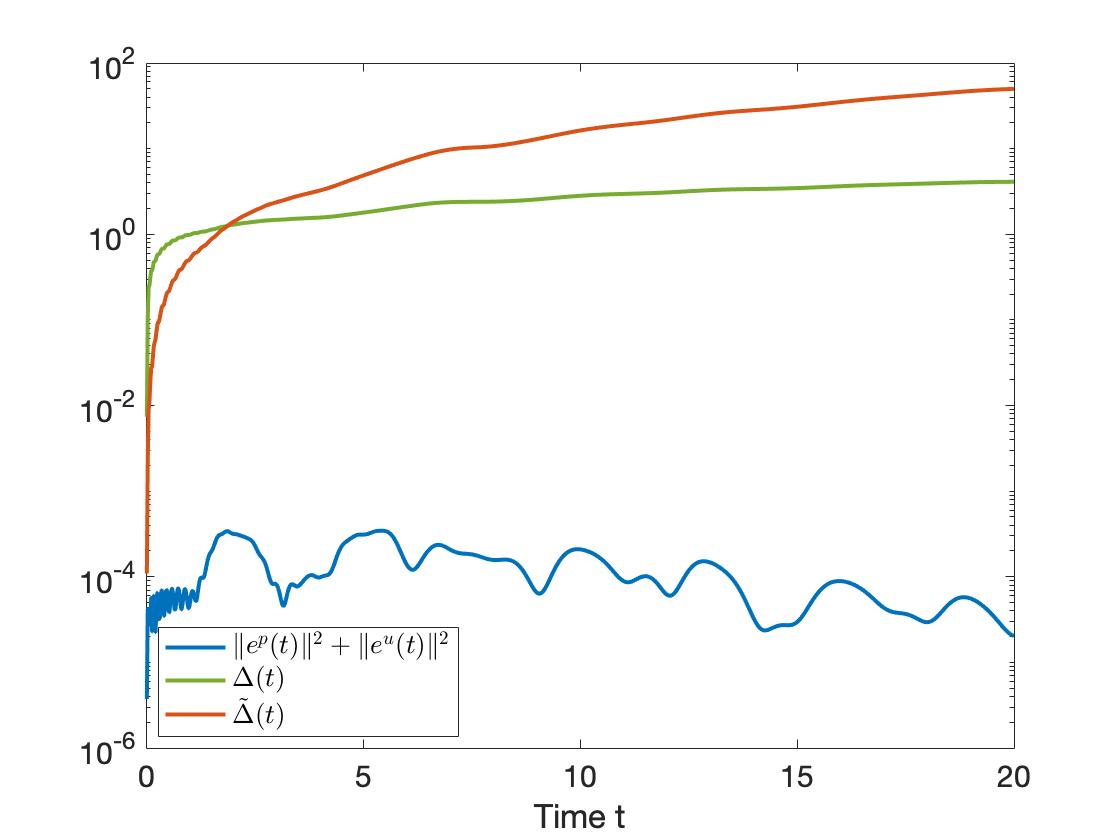}
\caption{Temporal evolution of actual error and bounds for $N = 53$ and $\mu = 2.3$.}\label{fig:samples}
\end{figure}

Figure~\ref{fig:samples} illustrates the temporal evolution of the actual error $\|e^p(t)\|^2 + \|e^u(t)\|^2$ and both error bounds $\Delta(t)$, $\tilde \Delta(t)$, $t\in[0,20]$ for the example of an intermediate friction coefficient $\mu=2.3$ and a reduced space dimension $N=23$. The results are qualitatively representative for the whole test sample and different model size. While both error bounds receive similarly large principal contributions at the beginning, our bound $\Delta$ grows much slower in time (due to the exponentially decaying term). The maximal value is reached at the end time, where it is one order of magnitude smaller for $\Delta$ than for $\tilde \Delta$.

\begin{figure}[thpb]
\includegraphics[width = \columnwidth]{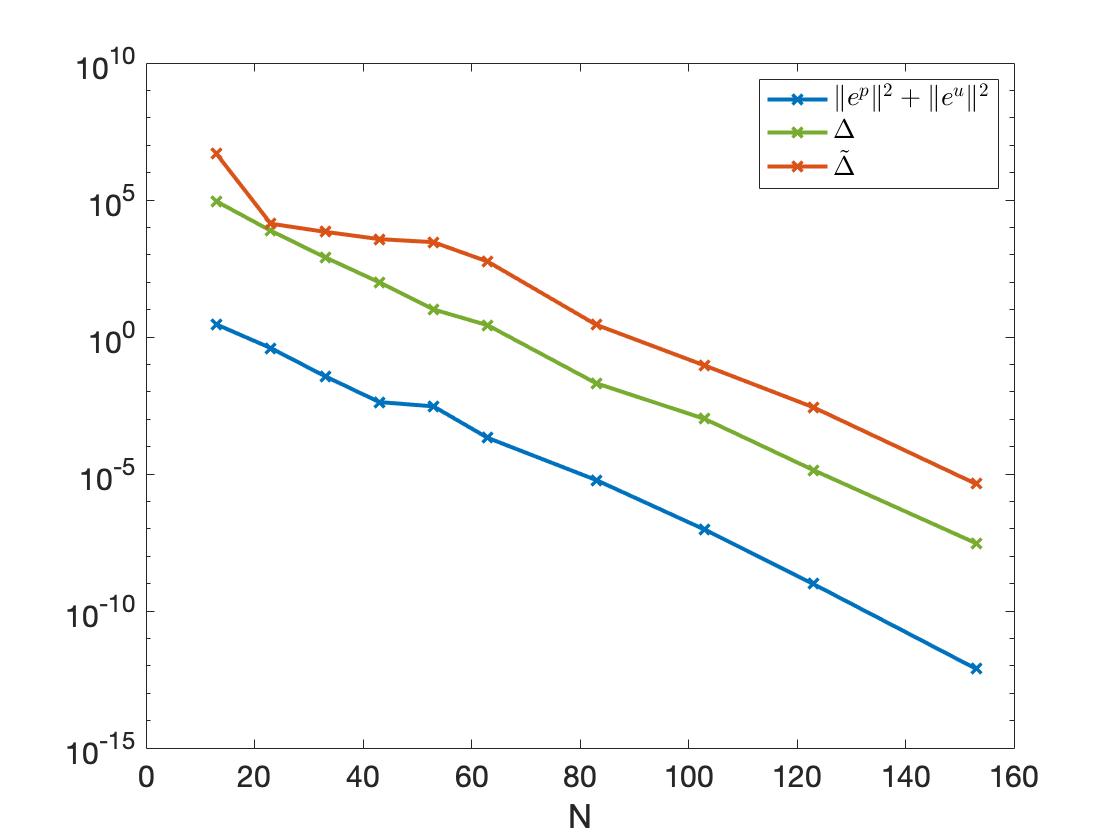}
\caption{Maximal values for error and bounds with respect to time and test sample for varying reduced space dimension.}\label{fig:bound_conv}
\end{figure}

\begin{figure}[thpb]
\includegraphics[width = \columnwidth]{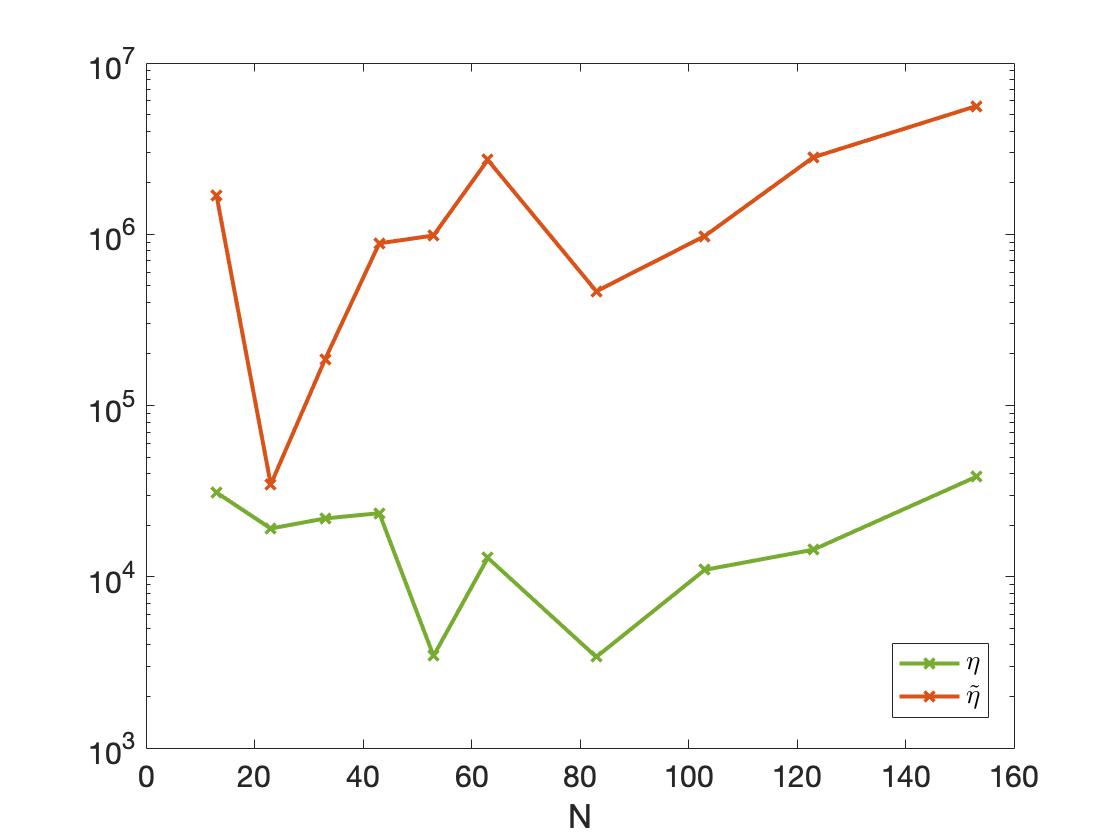}
\caption{Effectivities $\eta$ and $\tilde\eta$ for bounds $\Delta$ and $\tilde \Delta$. Maximal values with respect to time and test sample for varying reduced space dimension $N$ (cf., Fig.~\ref{fig:bound_conv}).}\label{fig:effective}
\end{figure}

Varying the model size by enriching the reduced spaces, the reduced approximations behave as desired, i.e., the error decreases for larger $N$. Figure~\ref{fig:bound_conv} shows for error and bounds the maximal values with respect to time and test sample for varying model size $N$. The qualitative behavior of both bounds is the same. Mimicking the behavior of the error, the bounds decrease linearly to zero for increasing~$N$. However, our bound $\Delta$ lies significantly below the one proposed by \cite{Haasdonk2011}, which makes it tighter. This property can also be concluded from the effectivity being the ratio of bound and error, i.e., 
\begin{equation*}
\eta(t) = \frac{\Delta(t)}{\|e^p(t)\|^2 + \|e^u(t)\|^2}.
\end{equation*}
(analogously $\tilde \eta$ for $\tilde \Delta$): the lower the effectivity, the sharper the bound. A perfect bound is indicated by an effectivity close to one. In Figure~\ref{fig:effective} the maximum effectivities with respect to time and test sample are visualized for varying model size. The maximal value for $\eta$ is in average two orders of magnitude smaller than the one for $\tilde{\eta}$. Thus, by exploiting the structure of the problem, our error bound $\Delta$ is tighter than the generic one $\tilde \Delta$.
For a desired approximation tolerance, it hence allows for lower-dimensional reduced models and faster convergence of the greedy selection algorithm (cf., Figure~\ref{fig:bound_conv}, where, e.g., a tolerance of $10^{-2}$ is ensured with $N=88$ for $\Delta$ and $N=115$ for $\tilde{\Delta}$). This makes our approach more efficient.

%%%%%%%%%%%
%%%%%%%%%%%
\section{Conclusion}

We derived an a posteriori error bound for a reduced basis approximation of the damped wave equations on networks. The resulting reduced order models were ensured to be stable, structure-preserving and efficient to evaluate. Our approach is based on the exponential decay of energy in the system, which allowed for sharp error bounds. This was verified numerically with comparison to other suitable error bounds. In future research, adaption of our error bound to optimization and data assimilation methods will be considered. 

\bibliography{stahl_mathmod_ref} 
\end{document}